\documentclass{amsart}
\usepackage[margin=1.0in]{geometry}
\usepackage{amsfonts,amsbsy,amsmath,amsthm,amssymb,verbatim}
\usepackage{url}
\usepackage{mathrsfs}
\usepackage{color}
\usepackage{comment} 

\usepackage[table]{xcolor}
\usepackage[normalem]{ulem}

\newtheorem{propo}{Proposition}[section]

\newtheorem{lemma}[propo]{Lemma}

\newtheorem{theo}[propo]{Theorem}

\newtheorem{remar}[propo]{Remark}

\newcommand{\F}{\mathbb{F}}
\renewcommand{\P}{\mathbb{P}}

\newcommand{\Q}{\mathbb{Q}}
\newcommand{\Co}{\mathbb{C}}

\newcommand{\bl}{\begin{lemma}\label}
\newcommand{\el}{\end{lemma}}

\newcommand{\ld}{,\ldots ,}
\newcommand{\ra}{ \rightarrow }

\newcommand{\GG}{\mathbf{G}}

\newcommand{\al}{\alpha}
\newcommand{\om}{\omega}
\newcommand{\si}{\sigma}

\newcommand{\dedication}[1]{%
  \begin{center}
    \large\itshape #1
  \end{center}
}

\newcommand{\bp}{\begin{proof} }
\newcommand{\enp}{\end{proof}}

\begin{document}
 
\title[]{Prime degree irreducible representations of simple algebraic groups 
and finite simple groups of Lie type}

\author{D.~L.~Flannery}
\address{University of Galway, Ireland}
\email{dane.flannery@universityofgalway.ie}
\author{A.~E.~Zalesski}
\address{Departamento de Mathem\'atica,
Universidade de Bras\'ilia, Bras\'ilia-DF, Brazil}
\email{alexandre.zalesski@gmail.com}  

\begin{abstract}
Let $r, p$ be primes and $k$ be a positive integer.  We show that if $G$ is a subgroup 
of $SL_r(p^k)$ lying only in the non-geometric Aschbacher class $\mathscr{C}_9$,  then 
$|G|$ is bounded above by a function of $r$ and $k$, independently of $p$. We also show 
that, up to conjugacy in $GL_r(p^k)$, the number of such $G$ is bounded above by a 
function of $r$ that is  independent of $p$ and $k$. Apart from being of interest in 
their own right, these results have an application in a computational version of the 
strong approximation theorem for finitely generated Zariski-dense subgroups of 
$SL_r(\mathbb{P})$, where $\mathbb{P}$ is a number field. 
\end{abstract}

\maketitle

\vspace{-5pt}

\dedication{Dedicated to the memory of Otto Kegel}

\vspace{-5pt}

\let\thefootnote\relax

\section{Introduction}

The study of finite linear groups is a central theme in group theory, going back 
to the origins of the subject. The following result of C.~Jordan (1870) 
is fundamental.

\smallskip 

{\bf Theorem\hspace{4pt}A} (\cite[Theorem~1.1, p.~444]{Feit})\textbf{.}
\emph{There exists a function $f(n)$ such that 
each finite subgroup $G$ of $GL_n(\mathbb{C})$ has an abelian normal subgroup whose 
index in $G$ does not exceed $f(n)$. }

\smallskip   

The optimal $f(n)$ is determined in \cite{CollinsI}.
Jordan's theorem holds more generally for $G<GL_n(\F)$ with $\F$ an arbitrary field of 
characteristic $p$ 
 coprime to $|G|$. It is false otherwise,
since $GL_n(p^k)<GL_n(\F)$ for all $k\geq 1$ and algebraically closed $\F$.

The following analogue of Theorem~A in positive characteristic is due to Brauer and Feit.

\smallskip 

{\bf Theorem\hspace{4.5pt}B} (\cite[Theorem~1.2, p.~444]{Feit}; 
cf.~\cite[Theorem~A]{CollinsIII})\textbf{.} 
\emph{Let $\F$ be a field of characteristic $p>0$. There exists a function 
$f(n,m)$ such that each finite subgroup $G$ of $GL_n(\F)$ with 
Sylow $p$-subgroups of order at most $p^m$ has an abelian normal subgroup whose index 
in $G$ does not exceed $f(n,m)$.}

\smallskip 

The prime characteristic case was investigated further by Larsen and Pink, who proved 
the following.

\smallskip 

{\bf Theorem\hspace{4pt}C} (\cite[Theorem~0.3]{LP})\textbf{.}  
\emph{Let $\F$ be a field of characteristic $p>0$. There exists a function $f(n)$ such 
that $|G:Z(G)|\leq f(n)$ for each finite quasisimple subgroup $G$ of $GL_n(\F)$,  unless 
$G/Z(G)$ is a simple group of Lie type in the defining characteristic $p$.}

\smallskip 

Theorems~B and C can be proved using the classification of finite simple groups and the 
lower bounds for minimal degrees of non-trivial irreducible projective representations 
of such groups that are available in \cite[$\S \, 5.3$]{KL}; see especially 
\cite[Table~5.3.A]{KL}. One can also use this information to calculate a function $f(n)$ 
as in Theorem~C explicitly.

Our results complement Theorems~B and C. For the case of Theorem~C where $G$ is of Lie type 
in defining characteristic and the degree $n$ is prime, we obtain, under certain restrictions,
 an upper bound on $|G|$ that is {\it independent of $p$}. In a sense, this can be viewed as 
a new kind of finiteness condition on groups of Lie type. 

To describe our results, we recall an approach to the classification of linear groups over
 finite fields developed by Aschbacher~\cite{A84}. 
This work was directed toward the classification of maximal subgroups of finite 
classical groups. It divides finite linear groups into several `classes' that are mostly 
of a geometric nature. Let $V=\F_q^n$, where $\F_q$ denotes the field of $q$ elements.  
We define sets $\mathscr{C}_1, \ldots , \mathscr{C}_9$ 
of subgroups $G<GL(V)\cong GL_n(q)$ as follows.

\vspace{4pt}

$\mathscr{C}_1$: $G$ is reducible. 
\vspace{2pt}

$\mathscr{C}_2$: $G$ is imprimitive. 
\vspace{2pt}

$\mathscr{C}_3$: $G$ stabilizes an extension field structure on $V$.
\vspace{2pt}

$\mathscr{C}_4$: $G$ preserves a tensor decomposition $U\otimes W$ of $V$ where 
$\mathrm{dim}(U), \mathrm{dim}(W)>1$ and $G$ does not permute $U$ and $W$.
\vspace{2pt}

$\mathscr{C}_5$: $G$ is conjugate to a subgroup of $GL_n(q_1)Z$, where $Z$ 
denotes the scalars of $GL_n(q)$ and $q = q_1^a$ for some integer $a > 1$.
\vspace{2pt}

$\mathscr{C}_6$: $n$ is a $t$-power for some prime $t$, $G$ is primitive, and $G$ 
has a normal extraspecial absolutely irreducible symplectic-type $t$-subgroup.
\vspace{2pt}

$\mathscr{C}_7$: $G$ preserves a tensor decomposition $U_1\otimes\cdots\otimes U_m$ 
of $V$, where $m>1$ and the $U_i$ for $i=1\ld m$ have the same dimension, with $G$ 
permuting the $U_i$. 
\vspace{2pt}

$\mathscr{C}_8$: $G$ normalizes the stabilizer of a non-degenerate bilinear or 
unitary form on $V$.
\vspace{2pt}

$\mathscr{C}_9$: $G/Z(G)$ is almost simple and $G\not \in \mathscr{C}_i$ for 
$i=1,\ldots,8$. 

\vspace{4pt}

\noindent Aschbacher proved that each maximal subgroup of $GL_n(q)$ lies in at least one
$\mathscr{C}_i$. The classes $\mathscr{C}_i$ are then defined for other classical 
groups $X_n(q) < GL_n(q)$ such as $SL_n(q)$, with the understanding that the 
$\mathscr{C}_i$ contain only proper subgroups of $X_n(q)$. 

Our definitions of the $\mathscr{C}_i$ are slightly less precise than those given in 
\cite[$\S\S \, 2.1, 2.2$]{Brayetal} (where $\mathscr{C}_9$ corresponds to the class
 denoted by $\mathscr{S}$). We note a characterization of Aschbacher classes via matrix rings: 
 $G\in \mathscr{C}_i$ for $i\not \in \{ 6,8,9\}$ if and only if 
$G$ normalizes a non-scalar proper subring of $Mat_n(\F_q)$; see \cite{Zr}.

The definitions of the $\mathscr{C}_i$ imply that if $G< GL_n(q)$ 
lies in $\mathscr{C}_9$ then 

\vspace{5pt}

 \hspace{15pt} the minimal non-central  normal subgroup $N$ of $G$ is absolutely 
irreducible, quasisimple, 

 \hspace{15pt} does not preserve a non-degenerate classical form on $V$, and  
is not realizable over any      \hfill  $(*)$

\hspace{15pt} proper subfield of $\F_q$.

\vspace{5pt}

\noindent Note that $N$ is exactly the intersection $G^\infty$ of all 
terms in the derived series of $G$. Indeed, $G/N$ embeds in $\mathrm{Out}(N/Z(N))$, 
so is solvable by Schreier's theorem.

Class $\mathscr{C}_9$ might be regarded as `exotic'  or irregular by contrast with the
other $\mathscr{C}_i$. Groups in the latter classes arise naturally, according to how they act
on $V$. For this reason, groups in $\mathscr{C}_i$ for $1\leq i\leq 8$ are usually easier 
to study, as they admit reduction of problems to groups of smaller degrees or over smaller 
coefficient fields. A discussion of the importance of class $\mathscr{C}_9$, along with a 
variety of applications and computational and probabilistic aspects, may be found in \cite{RD}.

Let $k$ be a positive integer and  $r, p$ be primes. Our first major result is the following.
\begin{theo}
\label{t11a}   
There exist functions $f_1(r,k)$ and $f_2(r)$ such that
\begin{itemize}
\item[{\rm (i)}] if $G<SL_r(p^k)$ lies in $\mathscr{C}_9$, then $|G|\leq f_1(r,k);$
\item[{\rm (ii)}] up to $GL_r(p^k)$-conjugacy, the number of subgroups $G$ of $SL_r(p^k)$ in 
$\mathscr{C}_9$ is at most $f_2(r)$.  
\end{itemize}
\end{theo}

We emphasize that neither $f_1$ nor $f_2$ depends on $p$, and $f_2$ does not 
depend on $k$. However, $f_1$ cannot be made independent of $k$ (see Section~\ref{Final}).
Since the functions $f_1,f_2$ are independent of $p$, Theorem~\ref{t11a} may (potentially)
 have some applications in Galois theory.  
 
The following companion result to Theorem~\ref{t11a}~(ii) gives an explicit bound in 
the key special case. 

\begin{theo}
\label{t11b}  Up to $GL_r(p^k)$-conjugacy, the number of subgroups $G$ of $SL_r(p^k)$ 
lying in $\mathscr{C}_9$ that are quasisimple of Lie type in the defining 
characteristic $p$ does not exceed $(2(3r)^{1/2}+1)\cdot r^{(r^2+8)/2}$.
\end{theo}

In fact, for quasisimple $G<SL_r(p^k)$ of Lie type in characteristic $p$ satisfying $(*)$, 
we will prove that $p< \allowbreak r^{r^2/2}$. This result is crucial in our proofs of 
Theorems~\ref{t11a} and \ref{t11b}. It yields immediately a coarse bound 
$|G|\leq |SL_r(p_0^k)|$, where $p_0$ is the greatest prime less than $r^{r^2/2}$ 
(cf.~\cite{Li}). The remainder of the paper is largely devoted to establishing the 
bound in Theorem~\ref{t11b} on the number of conjugacy classes of such $G$.
Again, our bound is far from being sharp. We leave the task of finding better bounds 
to future work.

Since $SL_2(p^k)=Sp_2(p^k)$ preserves a  bilinear form, $(*)$ excludes $r=2$. 
However, for our purposes degree $2$ is of no interest, because the non-solvable 
subgroups of $SL_2(p^k)$ are well-known. In general, we are not concerned with primes 
$r$ for which all quasisimple irreducible subgroups of $SL_r(p^k)$ have been determined, 
as Theorems~\ref{t11a} and \ref{t11b} can be verified by inspection of lists of these groups
(which could also provide tight values for $f_1$ and $f_2$ in the relevant degrees; 
cf.~\cite{Brayetal,HM,Lu}).  

The hypothesis of prime degree is essential. If $r$ is composite then 
Theorem~\ref{t11a}~(i) may fail to hold, as shown by tables in \cite[Chapter~8]{Brayetal}. 
Restricting the degree is an effective measure to facilitate the solution of 
problems in linear group theory. Many classifications of linear groups 
under such restrictions have been compiled. To cite just a few examples: finite 
irreducible monomial linear groups of prime degree over $\Co$ are classified 
up to conjugacy in \cite{DZ,BFO}; simple linear groups of relatively small degrees 
are listed in \cite{HM,Lu}, and irreducible projective representations of finite 
quasisimple groups of prime-power degree over $\Co$ are classified in \cite{BBO,MZ}.  

Theorem~\ref{t11a} extends \cite[Lemma~3.1]{DFH}. The proof in \cite{DFH} is incorrect, 
since it relies on the classification given by \cite[Theorem~1.1]{MZ} of 
projective representations of finite quasisimple groups in prime-power degree over 
$\Co$. The context and motivation for this earlier result is the development of 
a computational version of the strong approximation theorem for finitely generated 
Zariski-dense subgroups $H$ of $SL_r(\P)$, where $\P$ is a number field; 
cf.~\cite[Window~9]{LubotzkySegal} and \cite{DFHSAT}. The associated algorithm decides
 when each $\mathscr{C}_i$ contains a congruence image of $H$ modulo some 
maximal ideal of a finitely generated subring $R \subset \P$ such that $H\leq SL_r(R)$. 
Consequently, the algorithm finds the set of all primes $p$
modulo which $H$ does not surject onto $SL_r(p^k)$, for relevant $k\leq |\P:\Q|$. 
The classes $\mathscr{C}_1, \ldots ,\mathscr{C}_8$ may be eliminated 
straightforwardly;  $\mathscr{C}_9$ requires separate treatment. Possible congruence images 
in $\mathscr{C}_9$ are eliminated by means of an upper bound on the orders of subgroups of 
$SL_r(p^k)$ in $\mathscr{C}_9$. Of course, this order bound must be independent of $p$; while 
dependence on $k$ is not an issue. The existence of such a bound is guaranteed by 
Theorem~\ref{t11a}~(i). 
 
Due to Theorem~C, we focus on quasisimple absolutely irreducible subgroups $S<SL_r(p^k)$ 
of Lie type in characteristic $p$.  By Steinberg's theorem~\cite[Theorem~43]{St},
$S< \GG\leq SL_r(\F)$ where $\GG$ is a simple algebraic group of the same Lie type as $S$ 
over an algebraically closed field $\mathbb{F}$ of characteristic $p$. Let $V$ be the 
$\GG$-module afforded by this representation of $\GG$. To prove Theorems~\ref{t11a} and 
\ref{t11b}, we link our problem with the theory of Weyl modules for simple algebraic groups 
(see, for instance, \cite[$\S \, 2.1$]{Hum}). We use the fact that for each  irreducible 
representation of a simple Lie algebra $L$ over $\Co$, and for every prime $p$, there exists 
an indecomposable $\GG$-module $W$ in characteristic $p$ of the universal simple algebraic 
group $\GG$ of the same Lie type as $L$. This $\GG$-module is referred to as a \emph{Weyl module}. 
The dimension of $W$ is equal to the dimension of the irreducible representation of $L$.
Furthermore, $V$ is a composition factor of $W$ whose highest weight coincides with that
 of $V$. If $W$ is irreducible then $\dim V=\dim W$, so $L$ has an 
irreducible representation of prime degree $r$. Such representations are determined in 
\cite{Ka}. By Theorem~\ref{ka1} below, we then conclude

\vspace{6pt}

\hspace{50pt}  $\GG\cong SL(V)$ or $SO(V)$, or $\GG$ is of type $A_1$,
or $\GG$ is of type $G_2$ with $r=7$. \hfill $(**)$

\vspace{6pt}

\noindent  
The condition $(*)$ rules out these possibilities.
So we can assume that $W$ is reducible, and then we have the following
(to be proved in Section~\ref{PrimeDegreeRep}).
\begin{propo}
\label{wy1} 
Let $\GG$ be a simple algebraic group in characteristic $p>0$,
 $V$ an irreducible $\GG$-module of dimension $d$ with $p$-restricted highest weight 
$\omega$, and $W$ a Weyl module for $\GG$ of highest weight $\omega$. If $W$ is 
reducible then $p<\allowbreak \dim W<d^{d^2/2}$.
\end{propo}

Hence, if $(**)$ does not hold, then Proposition~\ref{wy1} bounds $p$ in terms of $d$.

\medskip  

\noindent \emph{Notation}. Background on algebraic groups, Lie algebras,
and their representation theory is given in \cite{Hm,MaT,St}. 
For the root system of a simple Lie algebra $L$ of rank $n$ or a simple algebraic 
group $\GG$ of rank $n$, we denote by $\Phi$ the set of roots, by $\Phi^+$ the set 
of positive roots with respect to simple roots $\al_1, \ldots , \al_n$, and by 
$\omega_1, \ldots , \omega_n$ the fundamental weights of the root system. 
A weight is then an expression $\om=\sum_i a_i\om_i$ with every $a_i$ an integer. 
The notation may be simplified by writing $(a_1\ld a_n)$ in place of $\sum_i a_i\om_i$.  
If $a_1\ld \allowbreak a_n$ are non-negative then  $\om$ is \emph{dominant}, and 
if $0\leq  a_1\ld a_n<p$ then $\om$ is \emph{$p$-restricted}.  There is a bijection
between the set of irreducible $L$-modules (respectively, irreducible $\GG$-modules)
 $V$ and the set of dominant weights, 
with the image of $V$ being called the \emph{highest weight} of $V$. 
We often write $V=V_\om$ to mean that $\om$ is the highest weight of $V$.

\section{Lie algebras}

\begin{theo}
\label{ka1} 
Let $L$ be a simple Lie algebra over $\Co$. Suppose that $V$ is an irreducible 
$L$-module with prime dimension $r$. Then one of the following holds:
\begin{enumerate}
\item $r>2$ and $L\cong \mathfrak{sl}_r$,
of type $A_{r-1};$ 
\item $r>2$ and $L\cong \mathfrak{so}_r$,
of type $B_{(r-1)/2};$    
\item  $r=7$ and $L$ is of type $G_2;$  
\item  $L$ is of type $A_1$, any $r$. 
\end{enumerate}
In all cases except $(1)$,
 $L$ preserves a non-degenerate symmetric bilinear form on $V$. \end{theo}

\bp See \cite[Theorem~1.6 and $\S \, 1.7.7$]{Ka}. The additional claim is well-known.
\enp

\bl{bn1} Let $L$ be a simple Lie algebra of rank $n$ over $\Co$, 
and let $V$ be an irreducible $L$-module 
with highest weight  $\om=(a_1, \ldots , a_n)$.  
Then $\dim V\leq (c+1)^l$, where 
$c=\mathrm{max}\{a_1, \ldots , a_n\}$ and $l=l(L)$ is the number of positive
 roots of $L$.\el

\bp The Weyl dimension formula for an irreducible representation of 
a semisimple Lie algebra with root system $\Phi$~\cite[p.~139, Corollary]{Hm}
yields $\dim V= \Pi_{\al\in\Phi^+}(1+m_\al)$, where 
 $m_\al=\frac{(\omega,\al)}{(\rho,\al)}$, $\rho=\allowbreak (1\ld 1)$, and 
$(\cdot \, ,\cdot)$ is the symmetric bilinear form on the weight lattice of $L$. 

Let $\al=\sum_i b_i\al_i$ where the $b_i$ are non-negative integers, and put
 $t_i=(\omega_i,\al_i)$.  Since $(\om_i,\al_j)=0$ for $i\neq j$, 
we have $(\rho,\al)=\sum_i t_ib_i$ and
$(\omega,\al)=\sum_i a_ib_it_i\leq c\sum_i b_it_i$, whence $m_\al\leq c$.
The lemma is now clear.
\enp

\section{Prime degree representations of simple algebraic groups}
\label{PrimeDegreeRep}
 
Each finite simple group of Lie type is the central quotient of a quasisimple group $G$ 
such that $G = \allowbreak \GG^\si:= \allowbreak \{ g \in \GG \mid g^\si = g\}$ for a 
simply connected simple algebraic group $\GG$ in prime characteristic, say $p$, and a 
Steinberg endomorphism $\si$ of $\GG$. 
The prime $p$ is called the defining characteristic of $\GG$ and of $G$.
We refer to $\GG^\si$ as a \emph{quasisimple group of Lie type}. Although $p$ is not always 
determined by $S$ as an abstract group, for any specified $\GG,\sigma$ as above such that 
$S=\GG^\si/Z(\GG^\si)$, we say that $p$ is the defining characteristic of $S$.

Let $H$ be a finite quasisimple group with $H/Z(H)$ simple of Lie type in defining characteristic $p$. 
Then $H$ is quasisimple of Lie type, unless $p$ divides $|Z(H)|$. The latter cases are exceptional and 
listed in \cite[Table~5, p.~xvi]{Atlas}.
Since $r$ is prime, if $H<SL_r(p^k)$  is irreducible then it is absolutely irreducible; so 
the scalar group $Z(H)$ has order coprime to $p$. We may therefore assume that $G$ in 
Theorems~\ref{t11a} and \ref{t11b} is a quasisimple group of Lie type.

For each simple Lie algebra $L$ over $\Co$ and each prime $p$, there exists a simple 
algebraic group $\GG$ constructed from $L$. Consequently we have nine families of 
simple algebraic groups, named by the corresponding Lie algebras. There are four 
classical families $A_n$ $(n\ge 1)$, $B_n$ ($n\ge 2$), $C_n$ $(n\ge 2)$, 
$D_n$ $(n\ge 3)$, where $n$ is an integer; and five exceptional types, denoted by 
$E_6$, $E_7$, $E_8$, $F_4$, and $G_2$. The subscript in a name is the \emph{rank} of 
$\GG$. For each $\GG$ and each algebraically closed field $\F$, there exists a unique 
universal group of points $\GG(\F)$ of $\GG$ over $\F$. Its center  $Z(\GG(\F))$ is 
finite but not necessarily trivial. The representation theory of $\GG(\F)$ does not 
depend on the choice of algebraically closed field $\F$ of characteristic $p$. So we 
may assume that $\F$ is the algebraic closure of $\F_p$. We always take $\GG$ to be 
the universal simple algebraic group of the given type.
 
\bl{p21}  
{\rm \cite[Table~2, Theorems~4.4 and 5.1, and Appendices~A.49--A.53]{Lu}}
Let $n$ be the rank of $\GG$, and let $V$ be an irreducible $\GG$-module of 
dimension $d>1$. The following hold.
\begin{enumerate}
\item 
$n\leq d-1$, $(d-1)/2$, $d/2$, $d/2$, for $\GG$ of type 
$A_n$, $B_n$ $(p > 2)$, $C_n$, $D_n$, respectively; except for type $B_2$, 
where $2=n\leq d/2$. 
\item 
$d\ge 27$, $56$, $248$, $25$, $6$ for $\GG$ of type 
$E_6$, $E_7$, $E_8$, $F_4$, $G_2$, respectively.
\item 
If $\GG$ is of classical type and $V$ is not a twist of the natural $\GG$-module 
then as second minimal dimension bounds we have $d\geq (n^2+n)/2$,
 $2n^2+n$ $(n\geq 7)$, $2n^2-n-2$, $2n^2-n-2$ $(n\geq 8)$, for $\GG$ of type 
$A_n$, $B_n$ $(n > 2)$, $C_n$, $D_n$, respectively.
\end{enumerate}
\el
 
\begin{remar}
{\em If $\GG$ is of type $B_n$ with $p=2$ then $n\leq d/2$ in (1); so
$n\leq d/2$ for all $n\geq 2$ and all $p$.}
\end{remar}

\begin{theo}
\label{gu12}  
Let $f_p(d)$ be the number of inequivalent $p$-restricted irreducible 
representations of $\GG$ of dimension at most $d$. Then $f_p(d)< d^4$.  
More precisely,
\begin{enumerate}
\item {\rm \cite[Theorem~3.2]{GLT}} if $p=2$ then $f_p(d)\leq d;$  
\item {\rm \cite[Theorem~2.14]{GLT}} if $\GG$ is of type $A_n$ and $p>2$ then 
$f_p(d)< d^{4};$
\item {\rm \cite[Theorem~4.2]{GLT}} if $\GG$ is not of type $A_n$ and $p>2$ then 
$f_p(d)\leq d^{5/2}$.
\end{enumerate}
\end{theo}

Note that while $p$ in Theorem~\ref{gu12} is fixed, the bounds on $f_p(d)$ are 
valid for all $p$. 

Our next goal is to bound $p$. 
Let $V=V_\omega$ be an irreducible $\GG$-module of dimension $d$ 
with highest weight $\omega$. By general theory adduced earlier, $V$ 
is a composition factor of a Weyl module 
$W=W_\omega$ for $\GG$ with highest weight $\omega$.  This module $W_\omega$ 
is indecomposable,  $\omega$ occurs with multiplicity $1$ in $W_\omega$, and 
$W_\omega$ has the same dimension 
as the highest weight irreducible module $W'$ 
for the simple complex Lie algebra of the same 
Lie type as $\GG$; see \cite[$\S\, 2$, p.~7]{Hum} (where $\overline{V}_\omega$ is  
used in place of $W_\omega$). 
The numbers of distinct weights of $W$ and $W'$ coincide. 
\bl{om1}
If $W$ is reducible then $p< \dim W$.\el

\bp A result of Jantzen~\cite[Theorem~II]{Ja} states that every 
$\GG$-module $M$ with $\dim M\leq p$  is completely reducible. Since $W$ is 
indecomposable, the lemma follows. \enp

Next, we bound $\dim W$.
\bl{in1} 
Let $V=V_\omega$ be an irreducible $\GG$-module of dimension $d$ 
with $p$-restricted highest weight $\omega$. Let $W=W_\omega$ be the Weyl 
module with highest weight	$\omega$. Then $\dim W < d^{d^2/2}$.
\el

\bp If $n$ is the rank of $\GG$ and $\omega=(a_1, \ldots , a_n)$, then $a_i \le d-1$ 
for all $i$.  Indeed, for every simple root $\al_i$ of $\GG$ there exists a simple algebraic 
subgroup $\GG_i\le \GG$ of type $A_1$ (see \cite[Theorem~8.17~(f)]{MaT}).
By \cite[Proposition~16.3]{MaT}, the restriction $V|_{\GG_i}$ has a composition factor
 $V_i$ with highest weight  $a_i$. Since $a_i<p$, we have 
$\dim V_i=a_i+1$~\cite[Remark~4.5]{Lu}. Therefore $ a_i+1\leq \dim V=d$ for all $i$, 
as claimed. 
Set $c= \allowbreak \max\{a_1, \ldots , a_n\}$, so $c\leq d-1$. 
By Lemma~\ref{bn1}, $\dim W\leq (c+1)^l\leq d^l$, where $l=|\Phi^+|$.
 
By \cite[Tables~I-IX]{Bo},  $l=(n^2+n)/2$, $n^2$, $n^2$, $n^2-n$, $36$, $63$, 
$120$, $24$, $6$ for $\GG$ of type 
$A_n$, $B_n$, $C_n$, $D_n$, $E_6$, $E_7$, $E_8$, $F_4$, $G_2$ respectively. 
In particular, $l\leq n^2$ for $A_n$, $B_n$, $C_n$, $D_n$, $E_6$; while 
$l< 2n^2$ in the other cases. So $\dim W< d^{2n^2}$ uniformly for all $\GG$.

To bound $l$ in terms of $d$, we recall from Lemma~\ref{p21} that the minimal 
dimension of an irreducible $\GG$-module is
 $n+1$, $2n+1$, $2n$, $2n$, $27$, $56$, $248$, $25$, $6$, 
respectively, except for $\GG$ of type $B_2$, where the minimal dimension is $4$. 

If $\GG$ is of type $A_n$ then $n\leq d-1$ and $d^l=d^{(n^2+n)/2}\le d^{(d^2-d)/2}$. 
 
If $\GG=B_n$ or $C_n$ then $d\geq 2n$ and $d^l=d^{n^2}\leq d^{d^2/4}$. 
 
If $\GG=D_n$ then $d\geq 2n$ and $d^l=d^{n^2-n}\leq d^{(d^2-2d)/4}$. 
 
Thus $\dim W < d^{d^2/2}$ uniformly for classical types. For the exceptional groups,
 $l < 2n^2$ and $n\leq d/3$; so $\dim W<d^{2d^2/9}<d^{d^2/2}$. Hence this bound 
is valid for all $\GG$.
\enp

\begin{remar} 
\label{Remarkin1}
{\em If $\dim W>1$ is minimal (see Lemma~\ref{p21}), then the structure of 
$W$ is well-known. In most cases $W$ is irreducible, and in the other cases the non-trivial 
composition factors are not of prime degree. Therefore, we can assume that $\dim W$ is not 
less than the second minimal dimension (indicated in  items (2), (3) of Lemma~\ref{p21}). 
This allows us to reduce the bound in Lemma~\ref{in1}.
For example, if  $\GG$ is of type $A_n$ then $d\geq (n^2+n)/2=l$ so $\dim W\leq d^l\leq d^d$.
Similarly,  $\dim W$ does not exceed $d^l\leq d^{d/2}$,  $d^{d-4}$, $d^{(d+2)/2}$, for $\GG$ 
of type $B_n$ $(n\geq 7)$, $C_n$, $D_n$ $(n\geq 8)$, respectively. If $\GG$ is of type 
$B_n$, $2\leq n\leq 6$, then $d\geq 2^{n}$;  and if $\GG$ is of type $D_n$, $3\leq n\leq 7$, 
then $d\geq 2^{n-1}$ (these bounds arise from the spin or half-spin representations).
For the exceptional groups one can use \cite{ZL}.}
\end{remar}
 
Now we specialize to prime degree. Here note that, for $r\leq 11$, there are no proper 
quasisimple subgroups of $SL_r(p^k)$ of Lie type in characteristic $p$ satisfying $(*)$; this 
is evident from the tables in \cite[Chapter~8]{Brayetal}.
\begin{theo}
\label{vv1} 
Let $r$ be a prime and $k$ be a positive integer. If $p$ is a prime such 
that $SL_r(p^k)$ contains an irreducible quasisimple subgroup $G$ satisfying the conditions
\begin{itemize}
\item $G$ is a quasisimple group of Lie type in characteristic $p$,
\item $G\neq SL_r(p^i)$ for all $i$ dividing $k$, and
\item $G$ does not preserve a non-degenerate symmetric bilinear form,
\end{itemize}
then $p < r^{r^2/2}$.
\end{theo}

\bp 
We assume that $r>3$ (see above). As explained at the beginning of Section~\ref{PrimeDegreeRep},
 we view $G$ as the image of an irreducible representation $\phi:G_1\ra GL_r(p^k)$, where 
$G_1=\GG^\si$ for a Steinberg endomorphism $\si$ of $\GG$. Let $V$ be the underlying space of $\phi$.
Since it has prime degree, $\phi$ is absolutely irreducible. By Steinberg's 
theorem~\cite[Theorem~43]{St}, $\phi$ extends to $\GG$. Since $r$ is prime, $\phi$ is 
tensor-indecomposable,  so we can assume that the highest weight $\omega=(a_1, \ldots , a_n)$ 
of $\phi$ is $p$-restricted; see \cite[Theorem~41]{St}. Let $W_\omega$ be the Weyl module for 
$\GG$ with highest weight $\omega$.
 
If $\dim W_\omega=\dim V_\omega$ then $\dim W_\omega=r$, and hence
$r$ is the dimension of an irreducible representation $\tau_\mu$ of the simple 
Lie algebra $L$ whose type is the same as the Lie type of $\GG$, and 
$\mu=(a_1, \ldots , a_n)$. 
(So the weights of $\phi$ and $\tau$ as strings of integers coincide.) 
By Theorem~\ref{ka1},  $L$ (and hence $\GG$) is either of type $A_{r-1}$, 
$\omega\in\{ \omega_1,\omega_{r-1}\}$; or of type 
$B_{(r-1)/2}$ with $\omega=\omega_1$; or of type  $A_1$ with $r\leq p$
and $\omega=(r-1)\omega_1$; or of type $G_2$ with $r=7$. In the former case, 
$\GG=SL_r(\mathbb{F})$ where $\mathbb{F}$ is the algebraic closure of $\F_p$. 
But this is contrary to the assumption. In the remaining cases, since $r>2$, we
see that $\phi(\GG)$ preserves a non-degenerate symmetric bilinear form on 
$V_\omega$ for this $\omega$.

Suppose that $V\neq W$. Then $W$ is reducible. By Lemmas~\ref{om1} and \ref{in1}, 
$p<\dim W_\omega< r^{r^2/2}$.
\enp
 
In the next section, we move on to proving Theorem~\ref{t11b}.

\section{Finite groups of Lie type}

Let $\GG$ be a simple algebraic group of universal type in characteristic $p>0$.  
Recall that the Steinberg endomorphisms of $\GG$ are classified (up to an inner 
automorphism multiple) in terms of a field parameter $p^t$ and the order $e$ of a 
Dynkin diagram symmetry of $\GG$~\cite[Theorem~22.5]{MaT}. Thus, a particular 
group $\GG^\si$ is identified by a pair $p^t,e$, and customarily denoted by 
${}^eG(p^t)$~\cite[Table~22.1]{MaT}. The superscript $e$ is dropped if $e=1$. 

\bl{ex1} 
\begin{enumerate}
\item None of the groups ${}^2B_2(2^{2m+1}),m>1$, $ {}^2F_4(2^{2m+1})$, and $ {}^2F_4(2)'$ 
has a $2$-modular irreducible projective representation of prime degree.
\item $G= {}^2G_2(3^{2m+1})$ has a $3$-modular irreducible projective 
representation of prime degree $r$ only for $r=7$. Every such representation of 
$G$ of degree $7$ is orthogonal.
\end{enumerate}
\el

\bp
The claim for ${}^2F_4(2)'$ follows by inspection of the list of the $2$-modular
 irreducible representation degrees in \cite[p.~188]{JLPW}. Let $G$ be any of the 
other groups. By \cite[Theorem~43]{St}, it suffices to prove the lemma for the 
algebraic group $\GG$ such that $G=\GG^\si$ for some Steinberg endomorphism $\si$ 
of $\GG$.  If $\phi$ is an irreducible representation of $\GG$ of prime degree with 
highest weight $\omega$, then as before we can assume that $\phi$ is $p$-restricted; 
thus $\omega =\sum_{i=1}^n a_i\omega_i$ with $0\le a_i\le p-1$ for all $i$.

If $G = {}^2B_2(2^{2m+1})$ with $m>1$ or $G = {}^2F_4(2^{2m+1})$,
then $a_i < 2$, and if $G = {}^2G_2(3^{2m+1})$ then $a_i < 3$.

The irreducible representation degrees of $ \GG=F_4$ with $a_i<2$ and $ \GG=G_2$ with
 $a_i<3$  are listed in \cite{LuebeckOnline} and \cite[p.~413]{GS}. The claims here 
then follow by inspection.
 
If $\GG$ is of type $B_2$ then $\omega\in\{(0,0),(1,0),(0,1),(1,1)\}$. 
By \cite[Appendix~A.22]{Lu},  the possible degrees of $\phi$ are
$1$, $4$, $4$, $16$, respectively. This completes the proof.
\enp

\begin{theo}
\label{t22} 
Let $G={}^eG(p^t)$, $t\geq 1$, be a quasisimple group of Lie type. 
Let $\phi:G\ra \allowbreak GL_r(\F)$ be an irreducible representation, where $r$ is 
prime and $\F$ is an algebraically closed field of characteristic $p$. Suppose that 
\begin{enumerate}
\item $\phi(G)<SL_r(p^k),$
\item  $\phi(G)$ is not conjugate to a subgroup $SL_r(p^i)$ of $GL_r(p^k)$ 
 that arises from the subfield embedding $\F_{p^i}\hookrightarrow \F_{p^k}$, for all 
$i$ properly dividing $k$.
\end{enumerate}
 Then $k=t$ or $k=et$. Hence for fixed $e$, the parameter $t$ is uniquely determined 
by $k$. \end{theo}

\bp For any finite 
irreducible group $H< GL_m(\F)$
there exists a least positive integer $l$ such that $H$ is conjugate in $GL_m(\F)$ to 
a subgroup of $GL_m(p^l)$~\cite[Theorem~3.4B]{Dixon}. In our setting, $m=r$ and $l=k$. 

Let $V$ be the underlying space of $GL_r(p^k)$. Since $r$ is prime, $V$ is an absolutely 
irreducible tensor-indecomposable $\mathbb{F}_{p^k}G$-module. 
By \cite[Proposition~5.4.6]{KL}, if $G$ is non-twisted then $t$ divides $k$ and 
(2) implies  that $k=t$.

Suppose that $G$ is twisted. By Lemma~\ref{ex1}, we can ignore 
$G\in\{{}^2B_2(2^{2m+1}), {}^2F_4(2^{2m+1}),{}^2G_2(3^{2m+1})\}$. The other twisted $G$ 
are ${}^2A_n(p^t)$, ${}^2D_n(p^t), {}^3D_4(p^t)$, and ${}^2E_6(p^t)$. 
By \cite[Proposition~5.4.6~(ii)~(b) and Remark~5.4.7~(a)]{KL}, either $p^t=p^k$ or 
$p^{2t}=p^k$, except when $G={}^3D_4(q)$, in which case $p^{3t}=p^k$ (see also 
\cite[Lemma~8.5]{ga}). That is, $k = t$ or $et$. 
\enp

\begin{theo}
\label{bb1} 
Let $r, p$ be primes and $k$ be a positive integer. There are at most  
$(6\cdot (3r)^{1/2}+7)\cdot r^4$ non-conjugate quasisimple groups $G$ of Lie type 
in defining characteristic $p$ such that
\begin{enumerate}
\item $G$  is an irreducible subgroup of $SL_r(p^k),$ 
\item up to conjugacy, $G$ is not contained in a subfield group 
$SL_r(p^i)$ for $i$ properly dividing $k$.
\end{enumerate}
\end{theo}
 
\bp 
As in the proof of Theorem~\ref{vv1}, we assume that $r>3$. Also, as explained at the beginning 
of Section~\ref{PrimeDegreeRep}, we can assume that $G=\phi(\GG^\si)=\phi({}^eG(p^t))$ for an irreducible 
representation $\phi$, where $\GG$ is a simply connected simple algebraic group and $\si$ is a 
Steinberg endomorphism of $\GG$. Thus we can use Theorem~\ref{t22}. Recall that by 
\cite[Theorem~43]{St}, $\phi$ extends to a representation of $\GG$; so we can use the results of 
Section~\ref{PrimeDegreeRep}.

In Theorem~\ref{t22}, the group $\GG$, representation $\phi$, and $e$ are 
all fixed. It remains to control what happens when $\GG$, $\phi$, and $e$ vary.

By Lemma~\ref{p21}, the rank $n$ of $\GG$ does not exceed  $r$; 
so the number of $\GG$ of classical type does not exceed $4r$. This bound has 
been improved in Remark~\ref{Remarkin1}, as we can assume that $\dim\phi$ is not
the dimension of the minimal non-trivial $\GG$-module.
By our observations there, and Lemma~\ref{p21}~(2), it follows that 
$r\geq (n^2+n)/3$ for all $\GG$. Hence $n<(3r)^{1/2}$, and the number of classical 
types of $\GG$ does not exceed $4\cdot(3r)^{1/2}$.  We add ${}^2A_n(p^t)$, 
${}^2D_n(p^t)$, the five untwisted exceptional types, ${}^2E_6(p^t)$  and 
${}^3D_4(p^t)$, obtaining at most  $6\cdot(3r)^{1/2}+7$ groups $G$.
(The other twisted groups are irrelevant due to Lemma~\ref{ex1}. 
For $r=7$ we might add ${}^2G_2(3^t)$, $t$ odd. However, for this 
and many other small primes $r$, the irreducible representations of degree $r$ are 
known: one can check by inspection of the data recorded in \cite{Lu} that the bound
 holds.)

We can assume that the highest weight $\omega$ of $\phi$ is $p$-restricted. 
Indeed, if $\psi:\GG\rightarrow GL_n(\F)$ is another irreducible 
representation with highest weight $p^j\omega$, 
then $\phi(G)$ and $\psi(G)$ are conjugate subgroups of $GL_n(\F)$.
Then by Theorem~\ref{gu12}, the number of inequivalent 
irreducible representations $\tau$ of a simple algebraic group $\GG$ of 
degree $r$ in characteristic $p$ does not exceed $r^4$. Note that the number of 
non-conjugate groups $\tau(\GG)$ does not exceed the number of 
inequivalent irreducible representations $\tau$.

By Theorem~\ref{t22}, for each simple algebraic group $\GG$ in 
characteristic $p$ and fixed $e$ there is at most one field parameter $t$ such 
that ${}^e G(p^t)$ has an irreducible representation of degree $r$ over $\F_{p^k}$ 
satisfying $(1)$ and $(2)$. Our total count is thus $(6\cdot (3r)^{1/2}+7).r^4$.
\enp

\begin{theo}
\label{th8}  
Let $r, p$ be  primes. For each integer $k\geq 1$, up to conjugacy in $GL_r(p^k)$ 
there are at most $(6\cdot (3r)^{1/2}+7)\cdot  r^{\frac{r^2}{2}+4}$ finite quasisimple 
groups $G$ of Lie type in defining characteristic $p$ such that $G<SL_r(p^k)$ and
 $G$ is not conjugate to a subgroup of $SL_r(p^i)$ for any $i$ properly dividing $k$.
\end{theo}

\bp 
As before we assume that $r>3$. By Theorem~\ref{vv1}, $p<r^{r^2/2}$. The number of 
such primes $p$ is about $(2r^{\frac{r^2}{2}-2})\big/ \log r$, but we choose instead 
the rough bound $r^{r^2/2}$. By Theorem~\ref{bb1}, for each $p$ there are at most 
$(6\cdot (3r)^{1/2}+7)\cdot r^4$ quasisimple groups $G$ of Lie type in defining 
characteristic $p$ that satisfy the conditions of the theorem. The result follows.
\enp

\section{Classical forms}

Recall that a Brauer character of a finite group is called \emph{real} if 
its values are real numbers.

\bl{sd1} If the Brauer character of an irreducible representation 
$\phi$ of a finite group $G$ is real, then $\phi(G)$ is contained 
in a symplectic or orthogonal group. 
\el

\bp The proof of \cite[Theorem~11.1, Chapter~IV, p.~189]{Feit} shows that $\phi(G)$ preserves 
a non-degenerate symplectic or skew-symmetric bilinear form on 
the underlying space of $\phi$. (Formally, \cite[Theorem~11.1, Chapter~IV, p.~189]{Feit} deals with 
characteristic $2$, but the reasoning is valid for arbitrary fields. In the notation of 
the proof of that result,  $c=\pm 1$ and $M'=\pm M$, so that $M$ is a Gram matrix of a 
symmetric or skew-symmetric form.) 
\enp

\bl{sd2} Let $G$ be a quasisimple group of Lie type in defining 
characteristic $p$,  and let $\phi:G\ra GL_m(p^t)$ be an irreducible representation  
of $G$. Suppose that $G$ is not of Lie type $A_n$ or $D_{2n+1}$ for $n>1$, nor $E_6$. 
Then $\phi(G)$ is contained  in a symplectic or orthogonal subgroup of 
$GL_m(p^t)$. 
\el

\bp By \cite[Proposition~3.1~(ii)]{TZ05}, each $p'$-element of $G$ is real,
i.e., conjugate to its inverse. Hence the Brauer character of $\phi$ is real, and
 the result follows from Lemma~\ref{sd1}. \enp

In Lemmas~\ref{sd1} and \ref{sd2}, we are referring to the full symplectic and 
orthogonal groups $Sp(V)$ and $O(V)$, where $V$ is the underlying space of $\phi$. 
In Lemma~\ref{sd2}, if  $m$ is a prime then  $\phi(G)\leq Sp(V)$ only for $m=2$.

\bl{ud2} Let $G$ be one of the groups ${}^2A_n(q)$ for $n>1$, ${}^2D_{2n+1}(q)$ for $n>2$, 
or ${}^2E_{6}(q)$, and let $\phi:G\rightarrow GL_m(p^k)$ be a tensor-indecomposable  
absolutely irreducible representation. Then $\phi(G)$ is contained in a proper 
classical subgroup of $GL_m(p^k)$. 
\el

\bp We can assume that $k$ is the least positive integer such that $\phi(G)$ is 
conjugate to a subgroup $GL_m(p^k)$. Let $G=G(q)$, $q=p^t$, and $V$ be the underlying 
space of $GL_m(p^k)$.
 As shown in the proof of Theorem~\ref{t22}, we can further assume that $k=t$ or $2t$. 
Let $\GG$ be the simple algebraic group such that $G=\allowbreak \GG^\si$. 
By \cite[Theorem~41]{St}, there exists a $\GG$-module $M$ with $q$-restricted highest 
weight $\lambda$, say, such that $V=M|_G$. 

Let $f$ be the graph automorphism of $\GG$ arising from a symmetry of the Dynkin diagram 
corresponding to $\GG$~\cite[Theorem~11.12]{MaT}. Then $f^2=1$ 
and $f$ permutes the weights of $M$. Denote by $M^{f}$ the $f$-twist of $M$.  
By \cite[Proposition~5.4.2~(iii)]{KL}, the highest weight of $M^{f}$ is $\tau(\lambda)$.
If  $k=t$ then  $M^f|_G\cong V$, and if $k=2t$ then $M^f|_G$ is the Galois conjugate of 
$V=M|_G$ corresponding to the Galois automorphism of $\mathbb{F}_{q^2}$ of order $2$ over 
$\mathbb{F}_q$. Therefore, we may write $V^f=M^f|_G$. In addition, $M(\tau(\lambda))$ is 
the dual of $M=M(\lambda)$. (Indeed, the dual $M^*$ of $M$ is of highest weight $-w_0(\lambda)$, 
where $w_0$ is the longest element of the Weyl group of $\GG$~\cite[Proposition~5.4.3]{KL}. 
As stated prior to \cite[Proposition~5.4.3]{KL}, $w_0$ acts on the weights $\mu$ of $\GG$ by 
sending $\mu$ to $-\tau(\mu)$, so $-w_0(\lambda)=\tau(\lambda)$.) Since  $M^*|_G=(M|_G)^*$,
it follows that $V^f\cong V^*$. So $V=V^f\cong V^*$ if $k=t$, while $V^*$ 
is isomorphic to the Galois conjugate of $V$ if $k=2t$. By \cite[Lemma~2.10.15]{KL}, 
$\phi(G)$ preserves a non-degenerate form on $V$ that is bilinear if $k=t$ and unitary if 
$k=2t$. Hence $\phi(G)$ is contained in a proper classical subgroup of $GL(V)$, namely, the 
stabilizer of this form. \enp

\section{Proof of Theorems~\ref{t11a} and \ref{t11b}}

Let $r>2$ and $p$ be primes and $G$ be a subgroup of $SL_r(p^k)$ satisfying $(*)$ as stated in 
the Introduction, where $N=\allowbreak G^\infty$. Put $S=N/Z(N)$. The groups
 in Lemma~\ref{ex1} can be excluded as they have no irreducible representations  of degree $r$.
Since $N$ is absolutely irreducible and $r$ is prime, $|Z(N)|\in \{1,r\}$.  Suppose that 
 $|Z(N)|=\allowbreak r$. Then $r$ divides the order $s$ of the Schur multiplier of $S$.
But for $r>3$, \cite[Table~1, p.~viii, and Table~5, p.~xvi]{Atlas} show that $r$ and $s$ are 
coprime, except when $S\in \{PSL_n(q), PSU_n(q)\}$ with $n>3$. In the latter event, 
$s \hspace{1pt} | \hspace{1pt} n$ and then $r \hspace{1pt} | \hspace{1pt} n$.
The minimal degree of a projective representation of such $S$ is $n$, except when $S=PSL_3(2)$ 
with $s=3$~\cite[$\S \, 5.3$]{KL}. Hence $n=\allowbreak r$ if $r>3$. 
By \cite[$\S \, 5.3$]{KL}, we have $p\hspace{1pt} | \hspace{1pt} q$, so $N\cong SL_r(q)$ 
and $G\in \mathscr{C}_5$: a contradiction. Thus, by the remark preceding Theorem~\ref{vv1}, we 
assume that $r>3$ and $N=S$ is simple. 

Observe that $G/Z(G)$ embeds in $\mathrm{Aut}(N)$ via conjugation action of $G$ on $N$, 
because $C_G(N)=Z(G)$. Thus $|G| \leq r |\mathrm{Aut}(N)|\leq \allowbreak c r |N| \cdot \log |N|$ 
where $c$ is an absolute constant (e.g., see \cite[pp.~334--335]{LubotzkySegal}).

First suppose that $S$ is not of Lie type in characteristic $p$. By Theorem~C, $|G|$ is 
bounded above by a function of $r$ only. 
Now a finite group of order $m$ has at most $m$ inequivalent irreducible representations  
over any field $\F$. Since distinct $GL_r(\F)$-conjugacy classes of the irreducible 
subgroups of $GL_r(\F)$ isomorphic to $G$ give rise to inequivalent irreducible 
representations of $G$ of degree $r$, the number of such conjugacy classes is at most 
$|G|$. So Theorem~\ref{t11a} holds in this case.

Henceforth, we suppose that $N$ is of Lie type in characteristic $p$. Then
$|N|< \allowbreak |SL_r(r^{kr^2/2})|$ by Theorem~\ref{vv1}, whence
Theorem~\ref{t11a}~(i) follows. 

Let $M$ be the normalizer of $N$ in $SL_r(p^k)$. Of course $G\leq M$ and $N$ is 
characteristic in $G$. Thus $|M/N|$ bounds the number of 
$GL_r(p^k)$-conjugacy classes of groups $G<SL_r(p^k)$ such that $G^\infty = N$, for 
fixed $N$. We first prove Theorem~\ref{t11b}, thereby bounding the number of possible 
$N$ up to conjugacy. Subsequently we bound $|M/N|$ by a function of $r$ only.  
Together these two bounds will complete the proof of Theorem~\ref{t11a}.

By $(*)$,  $N$ does not  preserve a non-degenerate classical form.
Lemmas~\ref{sd2} and \ref{ud2} then imply that the possible Lie types for $N$ are 
$A_n(p^t)$ or $D_{2n+1}(p^t)$ with $n>1$, or $E_6(p^t)$. Thus, for fixed $p$, we can 
replace $6\cdot (3r)^{1/2}+7$ in Theorems~\ref{bb1} and \ref{th8} by $2\cdot (3r)^{1/2}+1$. 
Theorem~\ref{t11b} now follows from Theorem~\ref{th8}. 

We proceed to bound $|M/N|$ in terms of $r$ only.
The conjugation action by $M$ on its absolutely irreducible simple subgroup $N$ 
induces an embedding of $M/NZ(M)$ onto a subgroup, $T$ say, of $\mathrm{Out}(N)$. 
As is noted in \cite[$\S\, 3$]{Atlas}, $\mathrm{Out}(N)$ is generated by three 
subgroups comprising so-called $d$-, $f$-, and $g$-automorphisms. Their orders are 
listed in the third, fourth, and fifth columns of \cite[Table~5, p.~xvi]{Atlas}. 
The $g$-automorphism group has order at most $6$. 
Denote the $d$- and $f$-automorphism groups by $D, F$, respectively. 
Then $D, F$ are cyclic, $DF=D\rtimes F$, and $DF\unlhd \mathrm{Out}(N)$.
 Furthermore, $|D|\leq 4$ unless $N\in \allowbreak \{PSL_n(q), PSU_n(q)\}$, 
in which cases $|D|\leq n$ and thus $n\leq r$, as $n$ is the minimal degree of 
a projective representation of such $N$. In view of Lemmas~\ref{ex1} and \ref{ud2}, 
$N$ is of non-twisted type. We can write $N=N(p^t)$, where $t=k$ by Theorem~\ref{t22}. 

Assume that $T\cap F=\{1\}$. Let $\al$ be any element of $T\cap DF$. Since 
$\langle \alpha\rangle$ acts by conjugation on $D$ and $|\mathrm{Aut}(D)| < r$,
we have $\alpha^{m} \in C_{DF}(D)= DC_F(D)$ for some $m$, $0<m<r$. Hence there 
exists $l$ such that $0<l\leq r$ and $\alpha^{ml} \in T \cap  F$; so $|\alpha|< r^2$.
Since $DF/D$ is cyclic, it follows that $|T\cap DF|<r^2 |D| \leq r^3$, and then 
$|M/NZ(M)| = |T|\leq 6 \cdot |T\cap DF| \leq 6r^3$. Therefore $|M/N| \leq  6r^4$.

To show that $T\cap F=\{1\}$, we suppose the contrary. Let $\phi:N\rightarrow N$ be the 
tautological  representation, and 
let $\al$ be a non-identity automorphism of $N$. Then $\al$ is realized by an element of 
$M$ if and only if $\phi$ is equivalent to the $\al$-twist $\phi^\al$ of 
$\phi$~\cite[Lemma~1.8.6]{Brayetal}. To see that this does not happen for any 
$f$-automorphism $\al$ of $N$, first extend $\phi$ to a polynomial representation 
$\overline{\phi}:\mathbf{N}\rightarrow GL_r(\F)$ of a simply connected simple
algebraic group $\mathbf{N}$, where $\F$ is an algebraically closed field of 
characteristic $p$.
Field automorphisms are restrictions to $N$ of the 
standard Frobenius endomorphisms $\sigma_i$ of $\mathbf{N}$ arising from the maps
$x\mapsto x^{p^i}$, $0\leq i \leq k-1$~\cite[p.~158]{St}. 
Since $\overline{\phi}$ is polynomial, for a field automorphism $\al=\sigma_i$ we have 
$\overline{\phi}(\sigma_i(g))=\sigma'_i(\overline{\phi}(g))$ for every $g\in \mathbf{N}$ 
and the $i$th standard Frobenius endomorphism $\sigma'_i$ of $GL_r(\F)$.  
The irreducible representations $\smash{\overline{\phi}}^{\hspace{1pt} \sigma_i}$,  
$0\leq i\leq k-1$, are inequivalent. Moreover,  since $N=N(p^k)$ is 
non-twisted and  $\phi$ is tensor-indecomposable, \cite[Theorem~43]{St}
implies that the restrictions of the $\smash{\overline{\phi}}^{\hspace{1pt} \sigma_i}$ 
to $N$ are inequivalent for distinct $i$; see \cite[p.~436]{Li}.
Hence $\phi$ and $\phi^\al$ are inequivalent, and so
 $T\cap F =\{1\}$ as claimed. We have now proved Theorem~\ref{t11a}~(ii).

\section{Concluding remarks}
\label{Final}

For the computational application, we need explicit bounds on the orders of \emph{all} 
$G < \allowbreak SL_r(p^k)$ in $\mathscr{C}_9$; i.e., for alternating 
$S= G^\infty/Z(G^\infty)$ and $S$ of Lie type in cross characteristic too.
These may be calculated using known facts (that hold in arbitrary degree). 
First, if $\F$ is any finite field and $\mathrm{Alt}(u)$ is a section of $GL_n(\F)$,
then $u\leq \allowbreak (3n+6)/2$~\cite[Proposition~10, p.~333]{LubotzkySegal}.
Secondly, \cite[Table~1]{SeitzZalesski} displays the least degree 
for which $S$ of Lie type in characteristic other than $p$ has a faithful projective 
representation over a field of characteristic $p$. These degree minima
bound the number of possible isomorphism types of $S$ independently of $p$.

We reiterate that the order bound function $f_1$ in Theorem~\ref{t11a}~(i) cannot be
 made independent of $k$. For a concrete example, let $\GG=A_2=SL_3$ and $p=7$. Then 
$\GG$ has an irreducible representation $\phi$ of degree $71$~\cite[Appendix~A.6]{Lu}. 
The highest weight of $\phi$ is $(2,5)$ or $(5,2)$, so $\phi$ and $\phi|_G$ are 
not self-dual. Thus, by \cite[Lemma~2.10.15]{KL},  if $G=SL_3(p^k)$ then
$\phi(G)$ does not preserve a  (non-degenerate) symmetric bilinear form, nor does it 
preserve a unitary form. Clearly $\phi(G)\cong SL_3(7^k)<SL_{71}(7^k)$ has
 unbounded order as $k\ra \infty$. However, the number of such groups up to 
$GL_{71}(7^k)$-conjugacy is bounded absolutely, by Theorem~\ref{t11b}. 

\bigskip

\bibliographystyle{amsplain}

\begin{thebibliography}{10}

\bibitem{A84}
M.~Aschbacher, On the maximal subgroups of the finite classical groups,
{\it Invent. Math.} 76 (1984), 469--514.

\bibitem{BBO} A. Balog, Ch.  Bessenrodt, J.  Olsson and K.  Ono, 
Prime power degree representations of the symmetric and alternating groups, 
{\it J. London Math. Soc.} 64 (2001), no.~2, 344--356.

\bibitem{BFO} Z. B\'acskai, D.~L.~Flannery and E.~A.~O'Brien,
 Classifying finite monomial linear groups of prime degree in characteristic zero,
Intern.  J.  Algebra and Comp.  31 (2021), 1547--1585.

\bibitem{Bo} N.~Bourbaki, {\it Groupes et alg\`{e}bres de Lie}, Ch. IV--VI, 
 Masson, Paris, 1981.

\bibitem{Brayetal}
J.~N.~Bray, D.~F.~Holt, and C.~M.~Roney-Dougal, \emph{The maximal subgroups of the 
low-dimensional finite classical groups}, London Math. Soc. Lecture Note Ser., vol.~407, 
Cambridge University Press, Cambridge, 2013.

\bibitem{CollinsI}
M.~J. Collins,  On Jordan's theorem for complex linear groups,
\textit{J. Group Theory} 10 (2007), no.~4, 411--423.

\bibitem{CollinsIII}
M.~J. Collins, Modular analogues of Jordan's theorem for finite linear groups,
\emph{J. Reine Angew. Math.} 624 (2008), 143--171.

\bibitem{Atlas} J.~H.~Conway, R.~T.~Curtis, S.~P.~Norton,  R.~A.~Parker, and 
R.~A.~Wilson, {\it An ATLAS of finite groups}, Clarendon Press, Oxford, 1985.

\bibitem{DFH}
A.~S. Detinko, D.~L. Flannery, and A.~Hulpke, Algorithms for experimenting with 
Zariski dense subgroups, {\it Exp. Math.} 29 (2020), 296--305.

\bibitem{DFHSAT}
A.~S. Detinko, D.~L. Flannery, and A.~Hulpke, The strong approximation theorem 
and computing with linear groups.
{\em J. Algebra} 529 (2019), 536--549.

\bibitem{Dixon}
J.~D. Dixon, {\em The structure of linear groups}, Van Nostrand Reinhold, London, 1971.

\bibitem{DZ}
J.~D. Dixon and A.~E. Zalesski, Finite imprimitive linear groups of prime degree, {\it J. Algebra}
276 (2004), 340--370.

\bibitem{Feit} W.~Feit, {\it Representation theory of finite groups}, North-Holland, 1982.

\bibitem{ga}  D.~Garzoni, Derangements in non-Frobenius groups, 
Proc. Lond. Math. Soc. (3) 131 (2025), no. 4, Paper No.~e70091, 31 pp. 

\bibitem{GS} D.~Gilkey and G.~Seitz, Some representations of exceptional 
Lie algebras, {\it Geom. Dedicata} 25 (1988), 407--416.

\bibitem{GLT} R.~Guralnick, M.~Larsen, and P.~H.~Tiep,  
Representation growth in positive characteristic and conjugacy classes
 of maximal subgroups, {\it Duke Math. J.} 161 (2012), no.~1, 107--137.
 
\bibitem{HM} G. Hiss and G. Malle, Low-dimensional representations of 
quasi-simple groups, {\it LMS J. Comput. Math.} 4 (2001), 22--63. 
Corrigenda: 5 (2002), 95--126. 

\bibitem{Hum} J.~E. Humphreys, {\em Ordinary and modular  representations 
of Chevalley groups}, Lecture Notes in Math. vol.~528, Springer-Verlag, New York, 1976. 

\bibitem{Hm} J.~E. Humphreys, {\em Introduction to Lie algebras and representation theory},
Springer-Verlag, New York, 1972. 

\bibitem{JLPW}  C.~Jansen, K.~Lux, R.~A.~Parker, and R.~A.~Wilson, {\it An ATLAS of
Brauer Characters}, Oxford University Press, Oxford, 1995.

\bibitem{Ja} J.~Jantzen, Low-dimensional representations of reductive groups 
are semisimple. In: {\em Algebraic groups and Lie groups},   {\em Austral.
  Math. Soc. Lect. Ser.}, vol.~9,  Cambridge University Press, Cambridge,
  1997, pp.~255--266.

\bibitem{Ka} N.~Katz, \emph{Exponential sums and differential equations}, 
Annals of Math. Studies, vol. 124, Princeton University Press, Princeton, New Jersey, 1999. 

\bibitem{KL} P.~B.~Kleidman and M.~W.~Liebeck, 
\emph{The subgroup structure of the finite classical groups},
 London Math. Soc. Lecture Note Ser., vol.~129, 
Cambridge University Press, Cambridge, 1990.

\bibitem{LP} M.~Larsen and R.~Pink, Finite subgroups of algebraic groups, 
\emph{J. Amer. Math. Soc.} 24 (2011), 1105--1158.

\bibitem{Li}  M.~W.~Liebeck,  On the orders of maximal subgroups of the finite classical
groups, {\it Proc. London Math. Soc.} (3) 50 (1985), no.~3, 426--446.

\bibitem{Lu}  F.~L\"ubeck, Small degree representations of finite Chevalley groups in 
defining characteristic,   {\it LMS J. Comp. Math.} 4 (2001), 135--169.

\bibitem{LuebeckOnline} F.~L\"ubeck, 
\url{https://www.math.rwth-aachen.de/~Frank.Luebeck/chev/WMSmall/F4-mod2.html}

\bibitem{LubotzkySegal}
A.~Lubotzky and D.~Segal, {\em Subgroup growth}, Birkh\"{a}user Verlag, Basel, 2003.

\bibitem{MaT}  G.~Malle and D.~Testerman, {\it Linear algebraic groups and 
finite groups  of Lie type}, Cambridge Stud. Adv. Math., vol. 133, Cambridge
  University Press, Cambridge, 2011.

\bibitem{MZ} G.~Malle and A.~E. Zalesski,
 Prime power degree representations of quasi-simple groups,
{\it Arch. Math.} (Basel) 77(2001), 461--468.

\bibitem{RD} C.~M.~Roney-Dougal, \emph{Maximal subgroups of  
  finite simple  groups: classification and applications}, 
	In {\em Surveys in Combinatorics} $2021$, London Math. Soc. Lecture Note Ser., vol.~470, 
Cambridge University Press, Cambridge, 2021, 343--369.

\bibitem{SeitzZalesski}
G.~Seitz and A.~E. Zalesski,
On the minimal degrees of projective representations of the finite
  Chevalley groups {II},
 {\it J. Algebra} 158 (1993), 233--243.

\bibitem{St} R.~Steinberg, {\it Lectures on Chevalley groups},  
Amer. Math. Soc. Univ. Lect. Series, vol. 66, Providence, Rhode Island, 2016.

\bibitem{TZ05} P.~H.~Tiep and  A.~E.~Zalesski,  Real conjugacy classes in 
algebraic groups and finite groups of Lie type, {\it J. Group Theory} 8 (2005), 
291--315.

\bibitem{Zr} Y.~Zarhin, Very simple representations: variations on a theme 
of Clifford, In: {\em Progress in Galois theory}, ed. H. Voelklein and T. Shaska, 
Springer, New York, 2005, 151--168.

\bibitem{ZL} Y.~Zhang and  Y.~Luan, Dimension formulas of the highest weight 
exceptional Lie algebra-modules, {\it AIMS Math.} 9, Issue 4, 
pp.~10010--10030.

\end{thebibliography}

\end{document}